\documentclass[review]{elsarticle}

\usepackage[dvipdfmx]{hyperref}
\usepackage{amsfonts}
\usepackage{amssymb,amstext}
\usepackage{amsthm}
\newtheorem{corollary}{Corollary}
\newtheorem{theorem}{Theorem}

\begin{document}

\begin{frontmatter}

\title{Concentration inequality using unconfirmed  knowledge}
\author[ntt]{Go Kato}
\address[ntt]{NTT Communication Science Laboratories, NTT Corporation, Atsugi-Shi, Kanagawa 243-0198, Japan}

\begin{abstract}
We give a concentration inequality 
based on the premise that random variables take values within a particular region.
The concentration inequality guarantees  that,
for any sequence of correlated random variables,
the difference between the sum of conditional expectations and that of the observed values takes a small value 
with high probability when the expected values are evaluated under the condition that the past values are known.
Our inequality  outperforms other well-known inequalities, e.g. the Azuma-Hoeffding inequality, especially in terms of the convergence speed when
the random variables are highly biased.
This high performance of our inequality is provided by the key idea in which we predict some parameters and 
adopt the predicted values in the inequality.
\end{abstract}


\end{frontmatter}

%\linenumbers

\section{Introduction}
For any correlated sequence of random variables, 
the sum of observed values is probably similar to that of the conditional expectations 
when the expected values are evaluated under the condition that
 the past observed values are known.
To express such a fact,
 we use an inequality that  gives the upper bound of the probability for atypical events, where 
the difference is bigger than a certain boundary value. We usually call such an inequality a concentration inequality.
This fact is equivalent to the fact that 
the difference between the $n$-th element and the first element of a martingale is small.  
We call such an inequality, which indicates the latter fact, a martingale concentration inequality.
Many relations\cite{A67,B04,DZ98,RS13} are known to be  martingale concentration inequalities.
However, we will use the former interpretation,
 because there are situations where
the conditional expectation itself has physical meaning~\cite{CXC14,DFR07,LNA19}, and the natural premises in the former case
are not always natural in the latter case.
Note that the condition on the expected values is necessary because we can easily define  correlated random values 
such that the sum of the observed values is far from that of the expected values with high probability. 

There are two types of relation between the sum of the conditional expectations and that of the observed values, i.e. concentration inequalities, that can be derived from the martingale concentration inequalities. However,  each type faces 
 a problem when we use them in applications.

The first type consists of  inequalities without statistical premises~\cite{A67,B04}. 
In this case,  we assume as a premises only that the random variables exist in a certain region, and we don't need to make assumptions about any statistical values, e.g. averages or variances.
One such martingale concentration inequality is the Azuma-Hoeffding inequality~\cite{A67}.
If all the observed values are given, we can directly check whether these premises are satisfied or not. Therefore, it is easy to guarantee that the premise holds. Although, such inequalities have  wide application,
  this type of concentration inequality has a defect.
To explain the defect, we will give some simple examples.
Assume that the random variables are almost independent of each other.
Accordingly, we can consider two situations.
The first one is that almost all the expected values are near the boundary of the region in which the random variables exist.
The second one is that  almost all the expected values are close to the median  of the region.
Therefore, the variances of the random variables must be small in the former  case, but
 may be large in the latter case.
As a result, 
we can see that there must be  a smaller upper bound of the probability of the atypical events in  the former case than in the latter case. 
However, since concentration inequalities without statistical premises depend only on the number of variables and the range of the region,
each concentration inequality, including the Azuma-Hoeffding inequality, give the same bound for the above two situations.
That means bounds given by such concentration inequalities for the former situation will be far from the achievable bound.

The second type consists of inequalities with statistical premises. 
In this case,  we assume as premise that certain statistical values, e.g. the average, the variance, etc., must be equal to or less than  certain constants~\cite{DZ98,RS13}.
In this case, even for the previous two situations, any inequality with statistical premises may give different bounds.
 Therefore, the bounds given by the inequalities may be good, and the problem explained in the previous paragraph 
 seems to be solved.
 However, a different problem appears in this case.
That is, in physical situations, it is usually difficult to  check the statistical values explicitly,
and 
 all we can do is  predict statistical values by random sampling. 
Therefore, in principle, we can treat only  observed values; We can never obtain exact statistical values.
That means,  concentration inequalities with statistical premises are inadequate for physical analyses except  in restricted situations.

In actual analyses, there are situations in which we know a priori that the expected values are  highly biased, as is the case of the former example above, though the bias  cannot be guaranteed~\cite{CXC14}. In such a situation, we could reluctantly use  the former type of inequality, rather than the latter type.

In this paper, we give a concentration inequality 
where we can use a priori knowledge under the premise about the region in which the random variables exist.
The inequality can be considered to be a generalization of the Azuma-Hoeffding inequality.
The inequality holds even if the knowledge is false.
If the  knowledge is true, the inequality gives a better bound.
For example, consider the case in which the random variables are almost independent and identically distributed, and  we can use the  average of the random variables as a priori knowledge. 
In this case, our upper bound is very similar to the achievable bound even if the average is near the boundary of the region.
This property may seem to be contradictory
in the sense that unconfirmed a priori knowledge can be used  to construct a confirmed better bound,
but we will show that this property is logically obtainable.

This paper is organized as follows.
In the next section, we formally express our main result and its corollaries.
In the third section, we show the key ideas to logically obtain the usefulness of unconfirmed knowledge.
In the fourth section, we give a proof of the main result.
In the fifth section, we compare our inequality and the known concentration inequality 
 in the case  where the random variables are independent.
The last section is devoted to a discussion.

\section{Main result}
Our main result is the following theorem:
\begin{theorem}
\label{th:main}
Let $\{X_m\}$ be a list of random variables, and 
$\{\mathcal F_m\}$ be a filtration that identifies random variables, including those variables $X_1\cdots X_m$; i.e. $\mathcal F_m$ is a $\sigma$ algebra which satisfies $\mathcal F_m\subseteq \mathcal F_{m+1}$ and $E(X_{m'}|\mathcal F_{m})=X_{m'}$ for $m'\leq m$.
Suppose that  the relation  $0\leq X_m\leq 1$ holds for any $m$.
In this case, for any $n\in \mathbb N$, $a\in\mathbb R$ and $b\in \mathbb R_{\geq 0}$,
\begin{eqnarray}
&&P\left(
\sum_{m=1}^n\left(E(X_m|\mathcal F_{m-1})-X_m\right)\geq
\left(b+a(2\frac{\sum_{m=1}^nX_m}n-1)\right)\sqrt n
\right)
\nonumber\\
&\leq&
\exp\left(-\frac{2(b^2-a^2)}{(1+\frac{4a}{3\sqrt n})^2}\right)
\label{eq:main_relation}
\end{eqnarray}
holds.
\end{theorem}

A characteristic of this inequality is that  the boundary value $\left(b+a(2\frac{\sum_{m=1}^nX_m}n-1)\right)\sqrt n $
depends on the random variables $\{X_m\}$. 
By rewriting the variables $\{X_m\}$ and the constants $a$ and $b$, the above relation gives other expression (\ref{eq:rel_01}) that enables us to 
imagine how to use this inequality efficiently.
Moreover, we can easily compare the rewritten expression with other concentration inequalities.

\begin{corollary}
\label{col:main}
Let  $\{X_m\}$ be a list of random variables, and
$\{\mathcal F_m\}$ be a filtration that identifies random variables, including those variables $X_1\cdots X_m$.
Suppose that
 $0 \leq X_m\leq 1$ holds for any $m$.
In this case,
for any  $n\in \mathbb N$, $\epsilon\in\mathbb R_{\geq0}$, $ \delta\in\mathbb R$  and $s\in\{1,-1\}$ such that $1-\delta(\delta-\frac{4s\epsilon}{3\sqrt n})> 0$, 
the following relation holds:
\begin{eqnarray}
&&P\left(
s(\sum_{m=1}^n(E(X_m|\mathcal F_{m-1})-X_m))\geq
\frac{1- \Delta(\delta-\frac{4s\epsilon}{3\sqrt n})}{1-\delta(\delta-\frac{4s\epsilon}{3\sqrt n})}
\epsilon \sqrt n 
\right)
\nonumber\\
&\leq&
\exp\left(-2\frac{\epsilon^2}{1-(\delta-\frac{4s\epsilon}{3\sqrt n})^2}
\right)
\label{eq:rel_01}
\end{eqnarray}
 where $\Delta:=1 -\frac{2}{n}\sum_{m=1}^n X_m$ is a random variable defined from $\{X_m\}$.
\end{corollary}
The relation for $s=1$ can be derived by substituting 
$\frac{\delta-\frac{4\epsilon}{3\sqrt n}}{1-\delta(\delta-\frac{4\epsilon}{3\sqrt n})}\epsilon $ and
 $\frac{1}{1-\delta(\delta-\frac{4\epsilon}{3\sqrt n})}\epsilon $ 
into $a$ and $b$, respectively.
The relation  for $s=-1$ is obtained from the one for $s=1$  by making the replacements
$X_m\rightarrow 1-X_m$, i.e. $\Delta \rightarrow - \Delta$, and  $\delta\rightarrow -\delta$.

If we set $\delta$ to $\frac{4s\epsilon}{3\sqrt n}$, the boundary value becomes a constant, and the above relation becomes just the concentration inequality derived from the Azuma-Hoeffding inequality.
To improve on this inequality, we can use a priori knowledge.  
 If the number of random variables $n$ is large enough, 
it is better to select the parameter $\delta $ such that $\Delta\simeq\delta$ with high probability.
For example, suppose we predict that the random variables $\{X_m\}$ are almost independent and the expected values of the random variables are almost equal to $p$. Then it is better to set the value $\delta$ as $1-2p$.
If the prediction is correct, we find that  
$
\frac{1- \Delta(\delta-\frac{4s\epsilon}{3\sqrt n})}{1-\delta(\delta-\frac{4s\epsilon}{3\sqrt n})}
$ is almost $1$ with high probability, and the obtained bound is probably much better than the one  given by the Azuma-Hoeffding inequality.

We should note  three things here. First, even if $\delta$ is far from the observed value of $\Delta $, the relation (\ref{eq:rel_01}) still holds. 
Second, if the random variables are independent, the conditional expectation $E(X_m|\mathcal F_{m-1})$
 is equal to the expected value $E(X_m)$. Therefore, the inequality indicates how small the difference between the sum of the expected values and that of the observed values is.
 Third,
 we can easily generalize the condition $0\leq X_m\leq 1$ to $\alpha_m\leq X_m\leq \beta_m$ in the case where all the values  $\beta_m-\alpha_m$ are a certain constant, which forces on us only a  trivial modification of the inequalities.

As mentioned  in the introduction, 
any martingale concentration inequality is equivalent to a corresponding concentration inequality.
Therefore, there is a martingale concentration inequality  equivalent to our inequality.
For convenience in making a comparison of our result with known martingale  inequalities, we will rewrite our result in the form of a martingale concentration inequality.

\begin{corollary}
\label{col:martingale}
Let $\{Y_m,\mathcal F_m\}_{m=0}^n$ be a real-valued martingale sequence, i.e. $E(Y_{m'}|\mathcal F_{m})=Y_{\min(m',m)}$.
Suppose that
there are  real random variables $\Delta_1,\cdots , \Delta_n$
which satisfy the relations 
$\frac 12(\Delta_m-1)\leq Y_m-Y_{m-1}\leq \frac 12(\Delta_m+1)$, and $E(\Delta_{m}|\mathcal F_{m-1})=\Delta_{m}$ for any $m\in\{1,2,\cdots ,n\}$. 
In this case,
for any  $n\in \mathbb N$, $\epsilon\in\mathbb R_{\geq0}$, $ \delta'\in\mathbb R$  and $s\in\{1,-1\}$ such that $1-\delta'(\delta'+\frac{2s\epsilon}{3\sqrt n})> 0$, 
the following relation holds:
\begin{eqnarray}
P\left(
s(Y_0-Y_n)
\geq
\frac{1-\Delta'(\delta'+\frac{2s\epsilon}{3\sqrt n})}{1-\delta'(\delta'+\frac{2s\epsilon}{3\sqrt n})}\epsilon\sqrt n 
\right)
&\leq&
\exp\left(-2\frac{\epsilon^{ 2}}{1-(\delta'+\frac{2s\epsilon}{3\sqrt n})^2}
\right)
%\nonumber\\
\label{eq:rel_11}
\end{eqnarray}
 where $\Delta':=
 \frac1n\sum_{m=1}^n  \Delta_m$.
\end{corollary}
In the case of $1\leq
\frac{\delta'+\frac{2s\epsilon}{3\sqrt n}}{1-\delta'(\delta'+\frac{2s\epsilon}{3\sqrt n})}\frac{2\epsilon}{\sqrt n} 
$ and 
$1>(\delta'+\frac{2s\epsilon}{3\sqrt n})^2$,
we can check that the left side of Eq. (\ref{eq:rel_11}) is  0  from the facts $-1\leq\Delta'\leq 1$ and $1-\delta'(\delta'+\frac{2s\epsilon}{3\sqrt n})> 0$.
In the case of $1\geq(\delta'+\frac{2s\epsilon}{3\sqrt n})^2$,
the right side of Eq. (\ref{eq:rel_11}) is trivially not less than 1.
In the other case,
 by substituting 
 $Y_m-Y_{m-1}+\frac12(1- \Delta_m)$ and
$\delta'+\frac{2s\epsilon}{\sqrt n}$
into  $X_m$ and $\delta$ respectively, the relation
(\ref{eq:rel_01})  becomes (\ref{eq:rel_11}).
Here, in order to apply Corollary \ref{col:main}, we have to assume the relation 
$1-
\frac{\delta'+\frac{2s\epsilon}{3\sqrt n}}{1-\delta'(\delta'+\frac{2s\epsilon}{3\sqrt n})}\frac{2\epsilon}{\sqrt n} 
>0$ since  $Y_0-Y_m+\frac{\epsilon\sqrt n(\delta'+\frac{2s\epsilon}{3\sqrt n})}{1-\delta'(\delta'+\frac{2s\epsilon}{3\sqrt n})} 
\Delta'$ is given as the value 
$\sum_{m=1}^n\left(E(X_m|\mathcal F_{m-1})-X_m\right)
+
\frac{\epsilon \sqrt n 
(\delta-\frac{4s\epsilon}{3\sqrt n})}{1-\delta(\delta-\frac{4s\epsilon}{3\sqrt n})}
 \Delta$
divided by $1-
\frac{\delta'+\frac{2s\epsilon}{3\sqrt n}}{1-\delta'(\delta'+\frac{2s\epsilon}{3\sqrt n})}\frac{2\epsilon}{\sqrt n} 
$.
Note that the random variables $\Delta_1,\cdots , \Delta_n$ identify the bias of the martingale.

 Therefore, if we can predict what $\frac 1n \sum_{k=1}^n  \Delta_k$ will be with high probability,
we can put that value into $\delta'$. If the prediction is correct, the inequality gives  with high probability a better bound than that of the Azuma-Hoeffding inequality.

\section{Key ideas}
To make unconfirmed information useful, we will exploit  the following three key ideas.

\subsection{The first key }
In the concentration inequalities derived from the conventional  martingale inequalities,
the boundary value, which identifies atypical events, is a constant depending only on the prior information.
However, this is not always necessary. 
In order to describe our first key idea,  we will consider a situation in which
the boundary value depends on the variances of the random variables.

With respect to the boundary value, the following two situations are equivalent to each other.
The first is, for a given boundary value, we evaluate the smaller upper bound of the atypical events' probability.
The second is,
for a given upper bound of the atypical events' probability, we evaluate the smaller boundary value.
In each situation, the bound of  the atypical events' probability becomes better.
As a result, we can change the boundary value depending on  the variances in order to give a good bound using a single inequality
in the case of small variances and in the case of large variances.

\subsection{The second key}
The first idea is not sufficient for our  purpose.
That is, usage of the variances contains the original problem 
affecting 
 concentration inequalities with statistical premises. That is, 
 statistical values like variance cannot be physically confirmed exactly.

To prevent such a problem, we can use the sum of the observed values rather than the variances.
If the sum of the observed values is very far from the median,
we can predict that the variances are very small.
Therefore, by using the sum of the observed values, we can implicitly use information about the variances
even if we can't guarantee what the variances are.

We will explicitly use a linear function of the sum of observed values as the boundary value.

\subsection{The third key}
The third key is the way to decide the linear function.
Depending on the choice of linear function, 
the upper bound of the atypical events' probability will become better. 
Therefore, we will use a priori knowledge to decide the function.
If the  knowledge is true, the bound becomes better.

Note that the knowledge is useful for finding  a better bound and 
that our inequality still holds even if the knowledge is false.

\section{Proof of Theorem \ref{th:main}}
In the case of $0\leq b\leq |a|$, the right side of eq. (\ref{eq:main_relation})  is not less than 1.
That means the relation (\ref{eq:main_relation}) trivially holds in such a case.
In the case of $a\leq -\frac {\sqrt n}{2}$ and $b>-a$, the left side of eq. (\ref{eq:main_relation}) 
is zero
from the assumption $0\leq X_m\leq 1$.
That means the relation (\ref{eq:main_relation}) also trivially holds in such a case.
 Therefore, we can use an additional  assumption $b> |a|$ and $a> -\frac {\sqrt n}{2}$ in this proof.

Since
$P(A\geq 0)\leq E(\exp(A\lambda ))$
for any non-negative value $\lambda$ and a random variable $A$, i.e. so-called Chernoff bound,
the left side of (\ref{eq:main_relation}) is upper bounded by the following values:
\begin{eqnarray}
&&
E\left(
\exp\left(\lambda(
-\beta+\sum_{m=1}^n(E(X_m|\mathcal F_{m-1})-\alpha X_m))
\right)
\right)
\nonumber\\
&=&
E\left(\left.
E\left(
\exp\left(\lambda(
-\beta+\sum_{m=1}^n(E(X_m|\mathcal F_{m-1})-\alpha X_m))
\right)
\right|\mathcal F_{n-1}\right)
\right)
\nonumber\\
&=&
E\left(
E\left(
\exp\left(\lambda
(
E(X_n|\mathcal F_{n-1})-\alpha X_n
)
\right)
|\mathcal F_{n-1}\right)
\right.
\nonumber\\
&&{}\times
\left.
\exp\left(\lambda(
-\beta+\sum_{m=1}^{n-1}
(E(X_m|\mathcal F_{m-1})-\alpha X_m))
\right)
\right)
\nonumber\\
&\leq&
\max_{0 \leq p\leq 1}
(p e^{\lambda (p-\alpha)}
+
(1-p)
 e^{\lambda p})
\nonumber\\
&&{}\times
E\left(
\exp\left(\lambda(
-\beta+\sum_{m=1}^{n-1}(E(X_m|\mathcal F_{m-1})-\alpha X_m))
\right)
\right)
\nonumber\\
&\leq&
\frac{1-e^{-\lambda\alpha}}\lambda
e^{\frac\lambda{1-e^{-\lambda \alpha}}-1}
E\left(
\exp\left(\lambda(
-\beta+\sum_{m=1}^{n-1}(E(X_m|\mathcal F_{m-1})-\alpha X_m))
\right)
\right)
\nonumber\\
&\cdots&
\nonumber\\
&\leq&
\left(\frac{1-e^{-\lambda\alpha}}\lambda\right)^n
e^{\frac{\lambda n}{1-e^{-\lambda \alpha}}-n-\beta \lambda}
\label{eq:rel_1}
\end{eqnarray}
where $\alpha:=1+\frac {2a}{\sqrt n}$ and $\beta:=\sqrt n(b-a)$.
In the first and second relations, we use the basic properties of conditional expectation:
 $E(E(A|\mathcal F_m))=E(A)$ and $E(A f(X_1,\cdots,X_m)|\mathcal F_m)=E(A|\mathcal F_m)f(X_1,\cdots,X_m)$
for any random variable $A$ and any function $f$, which comes from the property $E(X_m|\mathcal F_{m'})=X_m$ for $m\leq m'$.
In the third relation, i.e. the first inequality, we maximize the evaluated value with respect to 
the random variable $X_n$ under the constraint that $X_1\sim X_{n-1}$ are known and the average of $X_n$ is equal to $p$, i.e. $E(X_n|\mathcal F_{n-1})=p$, by
using the convexity of the function $e^{-x  \lambda\alpha}$, i.e. $\lambda\alpha\geq0$.
In the fourth relation, we evaluate the maximization for  $p\in\mathbb R$ without restricting region of $p$ by using the fact that $e^{-  \lambda\alpha}\leq 1$.
Note that, in the case of $\lambda=0$, we can evaluate $E(\exp(\lambda(
-\beta+\sum_{m=1}^n(E(X_m|\mathcal F_{m-1})-\alpha X_m))))$
 as $1$ and this is trivially upper bounded by the right-most expression $\alpha^n
e^{(\frac 1\alpha-1)n}=\lim_{\lambda\rightarrow 0}\left(\frac{1-e^{-\lambda\alpha}}\lambda\right)^n
e^{\frac{\lambda n}{1-e^{-\lambda \alpha}}-n-\beta \lambda}
$. 

The $\log$ on  the right side of eq. (\ref{eq:main_relation}) minus that of the above upper bound is 
\begin{eqnarray} 
&&-\frac{2(b^2-a^2)}{(1+\frac{4a}{3\sqrt n})^2}
+n+\beta \lambda-\frac{\lambda n}{1-e^{-\lambda \alpha}}
-n\log\left(\frac{1-e^{-\lambda\alpha}}\lambda\right)
\nonumber\\
&=&
n(-18\frac{x^2-y^2}{(4y+3)^2}
+1+\frac{(x-y)z}{2y+1}
\nonumber\\
&&{}
 -\frac{z }{(1-e^{-z})(2y+1)}
-\log\left(\frac{(1-e^{-z})(2y+1)}z\right))
\nonumber\\
&=:&nG(x,y,z)
\end{eqnarray}
where $x:=\frac b{\sqrt n}=\frac \beta n+\frac12(\alpha-1)$, $y:=\frac a{\sqrt n}=\frac12(\alpha-1)$, and $z:=\lambda \alpha$.
Therefore, all we have to do is check that, 
for any $x,y$ such that $x>|y|$ and $y>-\frac 12$, there is a non-negative value $z$ such that $G(x,y,z)\geq 0$. Note that $G(x,y,0)$ is defined to be $-18\frac{x^2-y^2}{(4y+3)^2}
+1 -\frac{1 }{2y+1}
-\log\left(2y+1\right)=\lim_{z\rightarrow 0}G(x,y,z)$, and 
$G(x,y,z)$ is an analytical function with respect to $x, z\in \mathbb R$ for $y>-\frac12$ even in the neighborhood of $z=0$.

In fact, this relation can be derived from the following four conditions:
\begin{eqnarray}
&& \forall y>-\frac12, \exists z \geq 0, \quad G(|y|,y,z)=0
\label{con:1}\\
&& \forall x\geq 0,\forall y, (\frac \partial{\partial x}G(x,y,z))|_{z=36 x\frac{2y+1}{(4y+3)^2}}= 0
\label{con:3}\\
&& \forall x\geq 0, \forall z\geq 0, \quad \frac {\partial^2}{\partial x^2}G(x,y,z)\leq 0
\label{con:4}\\
&& \forall x\geq 0,\forall y \;{\rm s.t.}\; -\frac12<y, \frac \partial{\partial x}G(x,y,36 x\frac{2y+1}{(4y+3)^2})\geq 0.
\label{con:5}
\end{eqnarray}
The above claim can be justified as follows.
From the first condition (\ref{con:1}), we can define a function $z_0(y)$ such that $G(|y|,y,z_0(y))=0$ for $y>-\frac12$.
Though $z_0(y)$ may be uniquely defined, it does not have to be in the following.
Here, what we have to do is use the uniquely defined function $z_0(y)$ for $y>-\frac12$.  
By using this function, we divide the region $\{(x,y)\;|\;|y|< x\;\land \;-\frac12<y\}$ into two regions:
\begin{eqnarray}
\makebox{region I}&&\{(x,y)||y|< x\;\land \; -\frac12< y\;\land \;x\leq 
\frac{(4y+3)^2}{36(2y+1)}z_0(y)
\}
\nonumber\\
\makebox{region II}&&\{(x,y)||y|< x\;\land \; -\frac12< y\;\land \;
\frac{(4y+3)^2}{36(2y+1)}z_0(y)
< x\},
%\nonumber\\
\end{eqnarray}
and we define the function $z_1(x,y)$ as:
\begin{eqnarray}
z_1(x,y)&:=&
\left\{
\begin{array}{cl}
z_0(y)& \makebox{in the case that $(x,y)$ is in region I}
\\
36 x\frac{2y+1}{(4y+3)^2}
& 
\makebox{in the case that $(x,y)$ is in region II.}
\end{array}
\right.
%\nonumber\\
\end{eqnarray}

When $(x,y)$ is in region I, we have 
\begin{eqnarray}
G(x,y,z_1(x,y))&=&G(x,y,z_0(y))
\nonumber\\
&\geq&G(|y|,y,z_0(y))=0.
\label{eq:region1}
\end{eqnarray}
The first and last equations are derived from the definitions of $z_1(x,y)$ and $z_0(y)$, respectively.
The second relation comes 
from 
$G(|y|,y,z_0(y))=0$,  ($\frac \partial{\partial x}G(x,y,z_0(y)))|_{x=
\frac{(4y+3)^2}{36(2y+1)}z_0(y)
}=0$, i.e. the second condition (\ref{con:3}), and the concavity of the function $G(x,y,z)$ with respect to $x$, i.e.  
 the third condition (\ref{con:4}).
When  $(x,y)$ is in region II, we have
\begin{eqnarray}
G(x,y,z_1(x,y))&=&G(x,y,
36 x\frac{2y+1}{(4y+3)^2}
)
\nonumber\\
&\geq&G(
\frac{(4y+3)^2}{36(2y+1)}z_0(y)
,y,z_0(y))
\nonumber\\
&\geq&G(|y|,y,z_0(y))=0.
\label{eq:region3}
\end{eqnarray}
The first and last equations are derived from the definitions of $z_1(x,y)$ and $z_0(y)$, respectively.
The second relation comes from the non-negative gradient of the function $G(x,y,
36 x\frac{2y+1}{(4y+3)^2}
)$ with respect to $x$, i.e.  the fourth condition (\ref{con:5}).
The third relation comes from 
$G(|y|,y,z_0(y))=0$,  ($\frac \partial{\partial x}G(x,y,z_0(y)))|_{x=
\frac{(4y+3)^2}{36(2y+1)}z_0(y)
}=0$, i.e. the second condition (\ref{con:3}), and the concavity of the function $G(x,y,z)$ with respect to $x$, i.e.  
 the third condition (\ref{con:4}).

What remains is to prove the four conditions (\ref{con:1})$\sim$(\ref{con:5}).

Two conditions (\ref{con:3}) and (\ref{con:4}) can be checked directly from the explicit expressions:
\begin{eqnarray}
\frac{\partial}{\partial x}G(x,y,z)&=&
-\frac{36x}{(4y+3)^2}
+\frac z{2y+1}\\
\frac{\partial^2}{\partial x^2}G(x,y,z)&=&
-\frac{36}{(4y+3)^2}
\end{eqnarray} 

Next, we check the first condition (\ref{con:1}).
In the case of  $y>0$: there is a value $ z_+> 0$ such that the relation $y=\frac { z_+}{2(1-e^{-z_+})}-\frac12$ holds.
We find 
\begin{eqnarray}
G(\frac { z_+}{2(1-e^{-z_+})}-\frac12,\frac { z_+}{2(1-e^{-z_+})}-\frac12,z_+)&=&0.
\end{eqnarray} 
Therefore,  the first condition holds.
In the case of  $-\frac12 <y<0$:  there is a value $ z_-> 0$ such that the relation $y=-\frac12+\frac { z_-}{2(e^{ z_-}-1)}$ holds.
We find 
\begin{eqnarray}
G(\frac12-\frac { z_-}{2(e^{ z_-}-1)},-\frac12+\frac { z_-}{2(e^{ z_-}-1)},z_-)&=&0.
\end{eqnarray} 
Therefore,  the first condition holds in this case as well.
In the case of  $y=0$: we can check 
\begin{eqnarray}
G(0,0,0)&=&0.
\end{eqnarray} 
Therefore,  the first condition holds in all three cases.

Finally, we check the last condition  (\ref{con:5}).
\begin{eqnarray}
&&\frac \partial{\partial x}G(x,y,
36 x\frac{2y+1}{(4y+3)^2}
)
\nonumber\\
&=&
((\frac 3w+w)\sinh w-3\cosh w)\frac{6(2y+1)}{(4y+3)^2\sinh w}
\nonumber\\&&{}+
(4w+2w\cosh 2w-3\sinh2w)\frac{3}{(4y+3)^2\sinh^2w}
+\frac{8w y^2}{(2y+1)(4y+3)^2}
\nonumber\\&&{}
\end{eqnarray}
where $w:=
18 x\frac{2y+1}{(4y+3)^2}
$.
We can check that all the terms are positive for $-\frac12<y$ and $w> 0$, i.e. $x>0$.
In particular, the positivity of the terms $(\frac 3w+w)\sinh w-3\cosh w$ and $4w+2w\cosh 2w-3\sinh2w$ for $w> 0$
can be shown as follows.
First, all the Taylor coefficients of the two terms at the point $w=0$ are positive,
\begin{eqnarray}
(\frac 3w+w)\sinh w-3\cosh w&=&\sum_{n=2}^\infty \frac{4n(n-1)}{(2n+1)!}w^{2n}
\nonumber\\
4w+2w\cosh 2w-3\sinh2w&=&\sum_{n=2}^\infty \frac{4^{n+1}(n-1)}{(2n+1)!}w^{2n+1}.
\end{eqnarray}
Second, these two terms are analytical functions of $w$ at any point in the complex plane; i.e. the Taylor series  converges for any $w$ case.
From the analyticity of the function $G(x,y,
36 x\frac{2y+1}{(4y+3)^2}
)
$ with respect to $x$ for $-\frac12<y$,  the non-negativity of the gradient is guaranteed even at $x=0$.

\section{Comparison}
We think that it would be informative to the reader to consider the case in which the random variables are independent of each other.

In statistical analyses,  we usually consider such cases and use the property that 
the sum of the observed values is similar to that of the expected values.
Actually, there are many inequalities that guarantee this property, such as the Hoeffding inequality~\cite{H63},
Bernstein's inequality~\cite{B24}, and Bennett's inequality~\cite{B62}.
Our inequality also guarantees this property in such situations.
That is, 
in the relations (\ref{eq:main_relation}) and (\ref{eq:rel_01}),
 $E(X_m|\mathcal F_{m-1})$ is equal to $E(X_m)$ 
when
the random variables are independent.
 We think it is still  unconventional that unconfirmed information, i.e. a priori knowledge, can be used to improve 
 the inequality.
In such situation, our inequality can be considered to be a generalization of the Hoeffding inequality.
As is the case with correlated random variables, if the sum of the expected values is far from the median, 
the Hoeffding inequality does not give a good bound. In contrast, our inequality can give a good one if we can use a priori knowledge about the expected values. Furthermore, if the a priori knowledge is true, our bound will be very similar to those derived from  
Bernstein's inequality or Bennett's inequality, which both explicitly use information about the variances of the random variables.

We can show another superiority when the average of the expected values is known.
In such a case, we usually  use the multiplicative Chernoff bound:
\begin{eqnarray}
P\left(s\sum_{m=1}^nX'_m <s(1-s\delta)np\right)<\left(\frac{e^{-s\delta}}{(1-s\delta)^{1-s\delta}}\right)^{np}
\end{eqnarray}
where $\{X'_m\}_m$  are independent random variables taking values in $[0,1]$, $p:=n^{-1}\sum_{m=1}^nE(X'_m)$, $0<\delta<1$ and $s\in\{1,-1\}$. By substituting $\frac\epsilon{\sqrt n p}$ into $\delta$, these relations can be rewritten in the equivalent form:
\begin{eqnarray}
P\left(s(np-\sum_{m=1}^nX'_m )>\epsilon \sqrt n\right)&<&
e^{-s\sqrt n \epsilon}\left(1-\frac{s\epsilon}{\sqrt n p}\right)^{-np +s\sqrt n \epsilon}
\nonumber\\
&=&e^{-\frac {\epsilon^2}{2p}+O(-n^{\frac12})}.
\end{eqnarray}
From Corollary \ref{col:martingale}, we can give a tighter bound than the above one as follows:
\begin{corollary}
\label{col:av}
Let  $\{X'_m\}$ be a list of independent random variables.
Suppose that
 $0 \leq X'_m\leq 1$ holds for any $m$.
In this case,
for any  $n\in \mathbb N$, $\epsilon\in\mathbb R_{\geq0}$ and $s\in\{1,-1\}$, 
the following  relation holds:
\begin{eqnarray}
P\left(s(np-\sum_{m=1}^nX'_m) \geq\epsilon \sqrt n\right)
&\leq&
\exp\left(-\frac{\epsilon^{ 2}}{2(p-\frac{s\epsilon}{3\sqrt n})({1-p+\frac{s\epsilon}{3\sqrt n}})}
\right)
\nonumber\\
&=&e^{-\frac {\epsilon^2}{2p(1-p)}+O(-n^{\frac12})}
\label{eq:col2}
\end{eqnarray}
 where $p:=\frac{1}{n}\sum_{m=1}^n E(X'_m)$.
\end{corollary}
In the case of $
\frac{\epsilon}{3\sqrt n}
\geq
\frac{2p(1-p)}{|1-2p|}
$,
the right side of Eq. (\ref{eq:col2}) is trivially not less than 1.
In the other case,
the relation (\ref{eq:col2})  is derived by substituting $\sum_{m'=1}^mX'_{m'}-E(X'_{m'})$, $1-2E(X'_m)$ and 
$1-2p$ into $
Y_m$, $\Delta_m$ and $\delta'$ respectively in Corollary \ref{col:martingale}.
Here, in order to apply Corollary \ref{col:martingale}, we have to assume the relation 
$\frac{\epsilon}{3\sqrt n}
<\frac{2p(1-p)}{|1-2p|}
$ since the relation is needed to guarantee the condition
$1-\delta'(\delta'+\frac{2s\epsilon}{3\sqrt n})> 0$.
Note that we don't have to use a priori knowledge to obtain superiority in this case.
It is fair to compare this corollary and the following similar relation derived form the recent result \cite{BK13,RS13}  (Lemma 2.2.3 in \cite{RS13}):
For any independent random variables $\{X'_m\}_m$ taking values in $[0,1]$ and any real positive value $\epsilon$,
\begin{eqnarray}
P\left(s\sum_{m=1}^n(E(X'_m)-X'_m) >\epsilon \sqrt n\right)
&<&
\exp\left(-\frac{\epsilon^2}{2\frac1n\sum_{m=1}^n c(E(X'_m),s)}\right)
%\nonumber\\
\end{eqnarray}
holds,
where 
\begin{eqnarray}
c(p,s)
&:=&
\left\{
\begin{array}{cl}
0&\makebox{if $p=0\lor p=1$}\\
\frac14&\makebox{if $p=\frac12$}\\
p(1-p)&\makebox{if $(0<p<\frac12\land s=1)\lor(\frac12\leq p\leq 1 \land  s=-1)$}\\
\frac{1-2p}{2\log{\frac{1-p}p}}&\makebox{if $(\frac12<p<1\land s=1)\lor(0<p<\frac12\land s=-1)$}.
\end{array}
\right.
\end{eqnarray}
We can easily check that Corollary \ref{col:av} has superiority in the case where $n$ is sufficiently large, 
all the expected values $E(X_m')$ are similar to the average of them and the average is  more (less) than $\frac12$ for $s=1(-1)$.

\section{Discussion and Conclusion}
We have given a concentration inequality for a series of random variables 
in the case where correlation may exist between the random variables and the range of the random variables is restricted to be from $0$ to $1$.
This inequality indicates that the sum of the observed values is similar to the sum of the conditional expectations.
This inequality enables us to use a priori knowledge about statistical values even if such  knowledge is uncertain.
As a result, by using a priori knowledge, we can get a desirable upper bound on the probability of atypical events.  
This bound is better than the well-known bound derived from the Azuma-Hoeffding inequality, especially in the case that
the sum of the observed values is far from the median.
Furthermore, the premise we use is the same one used in the Azuma-Hoeffding inequality.
Our inequality  is rather different   from other concentration inequalities in the sense that 
the boundary value, which identifies  atypical events, is defined as a function of random variables. 
 In order to select the function, a priori knowledge can be used even if 
there is no guarantee of its correctness. 
In fact, there are physical situations where we have some information about the random variables, though that information is not guaranteed to be true. We believe that  our inequality will be a powerful tool for analyzing such situations.

There remain two problems. 
We used a parameterized linear function
as a function of the boundary-value function. However, 
we think that there is a function that gives a better bound.
For example, if the upper bound of atypical events' probability is fixed, the envelope of the parameterized linear function is proportional to $\sqrt{1-\Delta^2}\sqrt n$ in the $n\rightarrow \infty$ limit.
This fact indicates that, if we use this function instead of a linear one, we may be able to get a better bound even
 without a priori knowledge.
The other problem is as follows.
In our inequality, 
we have specified the region in  which the random variables can exist.
 However,
in many concentration inequalities,
the  region may depend on the random variables.
If we could construct a concentration inequality for the case in which the  region may depend on the variables, it would have a very  wide range of application.

\section*{Acknowledge}
The author would like to thank  T. Kamaki and K. Azuma for valuable comments and encouragement.
The author also would like to thank G. C. Lorenzo for letting me conceive a key idea.
This work was partially supported by the JSPS Kakenhi (C) No. 17K05591.

\bibliographystyle{apsrev}

\end{document}